# Quantum-Assisted Stochastic Economic Dispatch for Renewables-Rich Power Systems

Xutao Han
*College of Electrical Engineering*
*Zhejiang University*
Hangzhou, China
xutao_han@zju.edu.cn

Zhiyi Li
*College of Electrical Engineering*
*Zhejiang University*
Hangzhou, China
zhiyi@zju.edu.cn

Yue Xu
*College of Electrical Engineering*
*Zhejiang University*
Hangzhou, China
yue_xu@zju.edu.cn

***Abstract*—Considering widely dispersed uncertain renewable energy sources (RESs), scenario-based stochastic optimization is an effective method for the economic dispatch of renewables-rich power systems. However, on classic computers, to simulate RES uncertainties with high accuracy, the massive scenario generation is very time-consuming, and the pertinent optimization problem is high-dimensional NP-hard mixed-integer programming. To this end, we design a quantum-assisted scheme to accelerate the stochastic optimization for power system economic dispatch without losing accuracy. We first propose the unified quantum amplitude estimation to characterize RES uncertainties, thereby generating massive scenarios by a few qubits to reduce state variables. Then, strong Benders cuts corresponding to some specific scenarios are selected to control the solution scale of Benders master problem in the iterative process, all of which are implemented by customized quantum approximation optimization algorithms. Finally, we perform numerical experiments on the modified IEEE 6-bus system to test the designed scheme.**

## I. INTRODUCTION

The penetration of renewable energy sources (RESs) (e.g., wind and solar energy) continues to soar to enable the low-carbon and economical operation of power systems. Thus, conventional power systems based mainly on fossil energy are being transformed into renewables-rich power systems, whose economic dispatch is more complicated with respect to uncertain RES generations [1]. To facilitate the economic dispatch of renewables-rich power systems, widely utilized scenario-based stochastic optimization can integrate uncertain RES generations to minimize operation costs [2]. However, since the generations of RESs dispersed in geographically distant regions are characterized as independent scenarios, the time and space complexities of scenario-based stochastic optimization will grow exponentially with the increasing number of RESs.

To generate massive scenarios, classical sampling methods based on the central processing unit (CPU) (e.g., sequential Monte Carlo and Latin Hypercube simulation) bear a burden of large computation to obtain accurate estimates, especially in the complicated distributions of RES generations. Also, the state space of stochastic optimization involving massive scenarios is very high-dimensional. Concurrently, just exploiting a small number of quantum bits (a.k.a., qubits) on quantum processing unit (QPU), quantum amplitude estimation (QAE) can generate an exponentially large number of scenarios to accurately characterize uncertainties [3]-[5]. The superiority of QAE has been proven in several other studies (e.g., financial risk analysis [3], numerical integration [4], and reliability analysis [5]) concerning diverse probability distributions. However, QAE has not yet been considered to accelerate the generation of requisite massive scenarios in stochastic optimization.

Besides, scenario-based stochastic optimization for the economic dispatch involves many binary variables (e.g., states of conventional generators and energy storages). CPU-based methods are almost impossible to directly solve such large-scale NP-hard mixed integer programming in limited iterations and time. Hence, in practical applications, scenario reduction (e.g., cluster) and binary variable relaxation are commonly pre-executed to reduce computing and running memory burdens [6]. However, both the two countermeasures may inevitably compromise decision accuracy and raise decision risks. Although the likelihood of obtaining inappropriate economic dispatch strategies is relatively low in a decision, the likelihood will continue to increase with the uninterrupted operation, concerning accumulative effects. Since any unsafe economic dispatch strategy may cause catastrophic consequences (e.g., transmission line overload and generator off-grid), an intuitive but effective countermeasure is to minimally reduce or even not reduce sampled scenarios and not relax binary variables. Fortunately, this contradictory dilemma on classical CPU can be resolved on QPU, leveraging the unique advantages of quantum approximation optimization algorithm (QAOA) to cope with binary variables. Nevertheless, QAOA is currently only utilized in conventional deterministic optimization (e.g., deterministic thermal unit commitment [7] and energy management [8]) and has not yet been considered for the stochastic optimization of economic dispatch in renewables-rich power systems.

To this end, we propose to assist the stochastic optimization by quantum computing, in which the massive uncertain scenarios can be efficiently generated by a few qubits and the pertinent NP-hard optimization process can be accelerated to polynomial time complexity without losing accuracy. The main contributions are summarized as:

(1) We design a quantum-assisted scheme for typical two-stage stochastic optimization under massive uncertain scenarios in the economic dispatch of renewables-rich power systems.

(2) We propose a unified QAE (UQAE) to characterize RES prediction errors and thus embed them in two-stage stochastic optimization with few qubits and lower estimation errors.

(3) We accelerate the multi-cut Benders decomposition that contains massive subproblems, in which the NP-hard feasible cut selection problem and master problem are reformulated in quantum-amenable forms and solved by the customized QAOA.

## II. UQAE-Based Uncertainty Characterization

### A. Modeling Distributions of RES Prediction Errors

The RES prediction errors (denoted as state variable $\xi$) generally follow normal or Gaussian mixture distributions according to different prediction algorithms [1]. As given in (1a), to model different error distributions, the continuous $\xi$ is discretized as $\xi(z_1)$ by $N_1+M_1$ qubits $z_{1,j}$, where $j$=−$M_1$,…,$N_1$; $M_1$ and $N_1$ are the number of qubits to represent the integer and decimal parts of $\xi$, respectively; $-(2^{N_1}-2^{-M_1-1})$ is the bias to cover both positive and negative errors.

$$\xi(z_1)=\sum_{j=-M_1}^{N_1} 2^j z_{1,j}-(2^{N_1}-2^{-M_1-1}) \tag{1a}$$

Then, the states of $|z_1\rangle=|z_{1,j}\rangle_{N_1+M_1}$ will follow the expected probability distribution after the designed unitary operator $\boldsymbol{P}$ in (1b), which is composed of several quantum gates (i.e., $RY(\theta)$ and CNOT) to simulate the probability $P(z_1)$. As shown in (1c), $RY(\theta)$ can rotate the Bloch vector of $z_{1,j}$ by angle $\theta$ around the Y-axis; CNOT can flip the target qubit once the control qubit is $|1\rangle$. In particular, Fig. 1 shows two intuitive examples to model the normal distribution (see Fig. 1(a)) and Gaussian mixture distribution (see Fig. 1(b)) by the above-mentioned qubits and quantum gates with different $\theta$ in $RY(\theta)$.

$$\boldsymbol{P}|0\rangle_{N_1+M_1}=\sum_{z_1=005\,0}^{115\,1}\sqrt{P(z_1)}|z_1\rangle \tag{1b}$$

$$RY(\theta)=\begin{bmatrix}\cos\frac{\theta}{2} & -\sin\frac{\theta}{2}\\ \sin\frac{\theta}{2} & \cos\frac{\theta}{2}\end{bmatrix},\ \text{CNOT}=\begin{bmatrix}1&0&0&0\\0&1&0&0\\0&0&0&1\\0&0&1&0\end{bmatrix} \tag{1c}$$

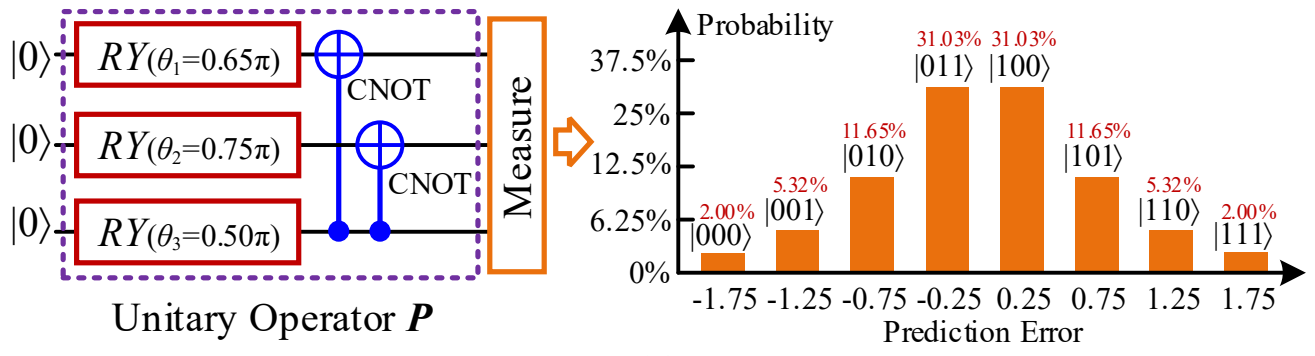

(a)

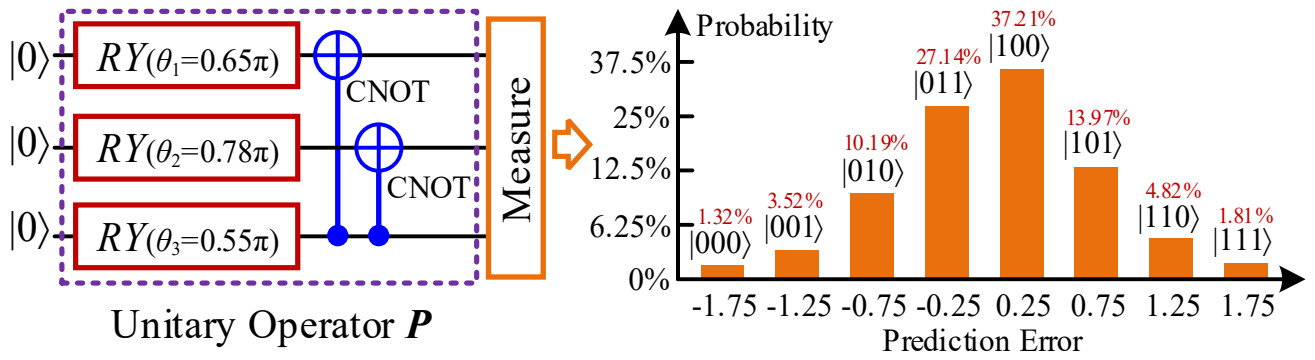

(b)

Fig. 1. Modeling different probability distributions of RES prediction errors. (a) Normal distribution; (b) Gaussian mixture distribution.

### B. Sampling RES Prediction Errors by UQAE

To implement UQAE, $\xi(z_1)$ in (1a) should be first unified as $\bar{\xi}(z_1)\in[0,1]$ in (2a). Then, the unitary operator $\boldsymbol{R}$ in (2b) acts on $|z_1\rangle$ and an ancilla qubit to discretize $\bar{\xi}(z_1)$.

$$\bar{\xi}(z_1)=\frac{\xi(z_1)+(2^{N_1}-2^{-M_1-1})}{2^{N_1+1}} \tag{2a}$$

$$\begin{aligned}|\psi\rangle &= \boldsymbol{R}(\boldsymbol{P}\otimes\mathbf{I})|0\rangle_{N_1+M_1}|0\rangle\\ &=\sum_{z_1=005\,0}^{115\,1}\sqrt{P(z_1)}|z_1\rangle(\sqrt{\bar{\xi}(z_1)}|1\rangle+\sqrt{1-\bar{\xi}(z_1)}|0\rangle)\end{aligned} \tag{2b}$$

where $\mathbf{I}$ is an identity operator; ⊗ represents tensor product.

To facilitate the design of a quantum circuit for UQAE, two orthonormal bases $|\psi_1\rangle$ and $|\psi_0\rangle$ are introduced to reformulate (2b) into (2c). Finally, unitary operators $\boldsymbol{Q}(2^0)$,…,$\boldsymbol{Q}(2^\alpha)$ in (2d) are exploited to sample $2^\alpha$ times with the gradually increasing probability to measure the state of ancilla qubit in $|1\rangle$. The sampling process for $\bar{\xi}(z_1)$ with probability $P(z_1)$ is equivalent to the amplitude estimation of $|\psi_1\rangle$ that integrates $\bar{\xi}(z_1)$ and $P(z_1)$. Fig. 2 shows the quantum circuit for UQAE including the unitary operators $\boldsymbol{P}$, $\boldsymbol{R}$, and $\boldsymbol{Q}$, which achieves quadratic speedup than classical sequential Monte Carlo simulation, whose maximum estimation gaps are $1/2^\alpha$ and $1/2^{\alpha-1}$, respectively.

$$|\psi\rangle=\sqrt{\sum_{z_1=005\,0}^{115\,1}P(z_1)\bar{\xi}(z_1)}|\psi_1\rangle+\sqrt{1-\sum_{z_1=005\,0}^{115\,1}P(z_1)\bar{\xi}(z_1)}|\psi_0\rangle \tag{2c}$$

$$\boldsymbol{Q}=(\mathbf{I}-2|\psi_0\rangle\langle\psi_0|)(\mathbf{I}-2|\psi\rangle\langle\psi|) \tag{2d}$$

in which,

$$|\psi_1\rangle=\frac{\sum_{z_1=005\,0}^{115\,1}\sqrt{P(z_1)}\sqrt{\bar{\xi}(z_1)}|z_1\rangle|1\rangle}{\sqrt{\sum_{z_1=005\,0}^{115\,1}P(z_1)\bar{\xi}(z_1)}} \tag{2e}$$

$$|\psi_0\rangle=\frac{\sum_{z_1=005\,0}^{115\,1}\sqrt{P(z_1)}\sqrt{1-\bar{\xi}(z_1)}|z_1\rangle|0\rangle}{\sqrt{1-\sum_{z_1=005\,0}^{115\,1}P(z_1)\bar{\xi}(z_1)}} \tag{2f}$$

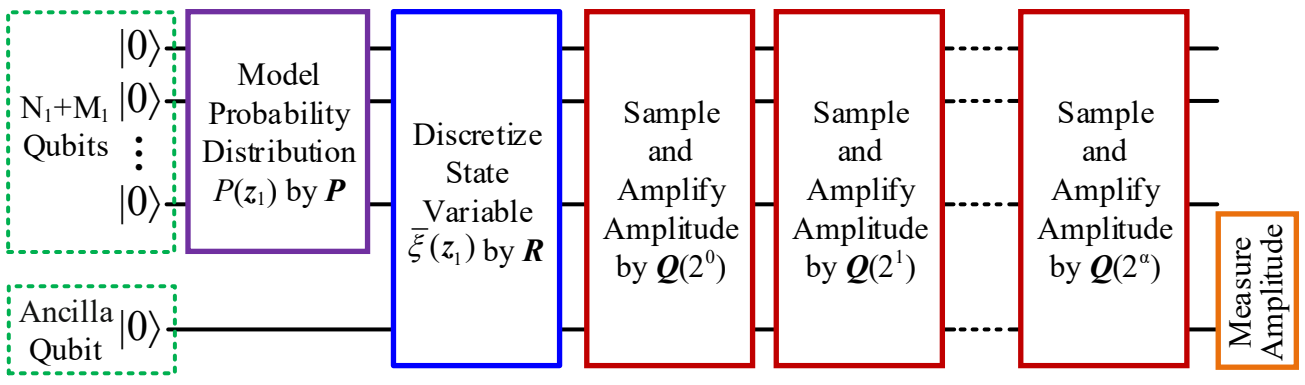

Fig. 2. Scheme of the quantum circuit for UQAE.

*Remark 1:* After measuring the amplitude of $|\psi_1\rangle$, $\bar{\xi}(z_1)$ is sampled and then easily transformed back to $\xi(z_1)$ according to (2a). Thus, QAE is generalized as UQAE to fit the case of uncertainty characterization with arbitrary value ranges of RES prediction errors, while the existing studies (e.g., [4] and [5]) are limited in applying QAE to estimate state variables in [0, 1].

## III. UQAE-QAOA Embedded Stochastic Optimization

### A. Quantum-Oriented Decomposition and Reformulation

Scenarios-based stochastic optimization for the economic dispatch of renewables-rich power systems is commonly formulated as a two-stage model (3). The first stage makes long-term decisions (e.g., the on/off states of conventional generators and the charging/discharging states of energy storages) based on predicted RES generations. After receiving the decisions from the first stage and embedding RES prediction errors by UQAE, fine-grained decisions are implemented to adapt current situation in the second stage (e.g., the outputs of conventional generators, the charging/ discharging power of energy storages, and thus the optimal power flow distribution). Since the second stage decisions are based on the trial solution from the first stage and will fed back to guide the correction of the first stage, the two stages are coupled with each other.

$$\min_{z_0, \boldsymbol{x}_s} \mathbf{a}^{\mathrm{T}} z_0 + \sum_{s=1}^{\mathrm{S}} \mathrm{p}_s \mathbf{b}_s^{\mathrm{T}} \boldsymbol{x}_s \tag{3a}$$

$$\text{s.t. } \mathbf{B} z_0 + \mathbf{C}_s \boldsymbol{x}_s \geq \mathbf{d}_s + \mathbf{A} \boldsymbol{\xi}_s \tag{3a}$$

$$z \in \{0,1\}^m, \boldsymbol{x}_s \in {}^{n}_{+}, \boldsymbol{\xi}_s = [\xi(z_1)]_{h\times 1} \leftarrow \text{UQAE} \tag{3c}$$

where $z_0$ is binary decision variables in the first stage; $\boldsymbol{x}_s$ and $\boldsymbol{\xi}_s$ are continuous decision and state variables under scenario $s$ in the second stage, respectively; $\mathrm{p}_s$ is the weight for scenario $s$ based on the samples of UQAE; S is the number of scenarios; $m$, $n$, and $h$ represent the dimension of variables; $\mathbf{a}$, $\mathbf{b}_s$, $\mathbf{d}_s$, $\mathbf{A}$, $\mathbf{B}$, and $\mathbf{C}_s$ are constant vector/matrix with appropriate dimensions.

Benefiting from provable convergence, multi-cut Benders decomposition is a prevailing algorithm exploited for solving model (3) by yielding a master problem and massive independent subproblems. The linear programming (LP) subproblems are suitable for solving on CPU but coupled with massive scenarios. The master problem is NP-hard mixed integer linear programming for CPU, whose time and space complexities rapidly increase with optimal and feasible cuts fed back from subproblems in iterations. Hence, we propose to reformulate the master problem and subproblems as quantum-oriented (4) and (5), respectively, and thus embed UQAE and QAOA modules to accelerate the accurate solution of model (3).

$$\min_{z_0, z_{2,i}} \mathbf{a}^{\mathrm{T}} z_0 + \sum_{i=-\mathrm{M}_2}^{\mathrm{N}_2} 2^i z_{2,i} \tag{4a}$$

$$\text{s.t. } \sum_{i=-\mathrm{M}_2}^{\mathrm{N}_2} 2^i z_{2,i} \geq \sum_{s=1}^{S_{\mathrm{op}}} \mathrm{p}_s (\mathbf{d}_s + \mathbf{A}\boldsymbol{\xi}_s - \mathbf{B} z_0)^{\mathrm{T}} \boldsymbol{u}_{s,\mathrm{op}}^* \tag{4b}$$

$$0 \geq (\mathbf{d}_s + \mathbf{A}\boldsymbol{\xi}_s - \mathbf{B} z_0)^{\mathrm{T}} \boldsymbol{u}_{s,\mathrm{fea}}^* \tag{4c}$$

$$z_0 \in \{0,1\}^m, \boldsymbol{\xi}_s = [\xi(z_1)]_{h\times 1} \leftarrow \text{UQAE}, z_{2,i} \in \{0,1\} \tag{4d}$$

$$\max_{\boldsymbol{u}_s} (\mathbf{d}_s + \mathbf{A}\boldsymbol{\xi}_s - \mathbf{B} z_0^*)^{\mathrm{T}} \boldsymbol{u}_s \tag{5a}$$

$$\text{s.t. } \mathbf{C}_s^{\mathrm{T}} \boldsymbol{u}_s \leq \mathbf{b}_s \tag{5b}$$

$$\boldsymbol{u}_s \in {}^{l}_{+}, \boldsymbol{\xi}_s = [\xi(z_1)]_{h\times 1} \leftarrow \text{UQAE} \tag{5c}$$

where $\mathrm{N}_2$+$\mathrm{M}_2$ qubits $z_{2,i}$ with the weight of $2^i$ relax the sum of optimal solutions for all subproblems; $z_0^*$ is the trial solution of master problem; $\boldsymbol{u}_{s,\mathrm{op}}^*$ and $\boldsymbol{u}_{s,\mathrm{fea}}^*$ are trial optimal and feasible solutions, respectively, about the dual variable $\boldsymbol{u}_s$ of subproblem $s$; $S_{\mathrm{op}}$ is the number of accumulated optimal cuts; $l$ is the dimension of $\boldsymbol{u}_s$.

*Remark 2:* Considering the specific characterization of two-stage stochastic optimization, the reformulation in (4) and (5) is different from existing deterministic optimization by Benders decomposition on QPU+CPU (e.g., [9] and [10]). In particular, only we have embedded RES prediction errors by UQAE, and thus massive subproblems (5) are formulated under diverse scenarios while existing studies only involve one subproblem.

Then, we transform the integer linear programming (ILP) (4) into quadratic unconstrained binary optimization (QUBO) (6) by Hamiltonian $H_{\mathrm{MP}}$, which is composed of three Hamiltonians $H_1$, $H_2$, and $H_3$ that correspond to objective (4a), aggregated optimal cut (4b), and independent $S_{\mathrm{fea}}$ feasible cuts (4c), respectively. After setting appropriate penalty factors $\lambda_1$ and $\lambda_2$ in (6c) and (6d), searching for the lowest energy of $H_{\mathrm{MP}}$ is equivalent to the optimization process of ILP (4) that minimizes the objective (4a) constrained by (4b)-(4d).

$$\min_{z_0, z_{2,i}, z_{3,r}, z_{4,ku}} H_{\mathrm{MP}} = H_1 + H_2 + H_3 \tag{6a}$$

in which,

$$H_1 = \mathbf{a}^{\mathrm{T}} z_0 + \sum_{i=-\mathrm{M}_2}^{\mathrm{N}_2} 2^i z_{2,i} \tag{6b}$$

$$H_2 = \lambda_1 \left( \begin{array}{c} \sum_{i=-\mathrm{M}_2}^{\mathrm{N}_2} 2^i z_{2,i} - \sum_{r=-\mathrm{M}_3}^{\mathrm{N}_3} 2^r z_{3,r} - \\ \sum_{s=1}^{S_{\mathrm{op}}} \mathrm{p}_s (\mathbf{d}_s + \mathbf{A}\boldsymbol{\xi}_s - \mathbf{B} z_0)^{\mathrm{T}} \boldsymbol{u}_{s,\mathrm{op}}^* \end{array} \right)^2 \tag{6c}$$

$$H_3 = \lambda_2 \sum_{k=1}^{S_{\mathrm{fea}}} \left( \sum_{u=-\mathrm{M}_4}^{\mathrm{N}_4} 2^u z_{4,ku} + (\mathbf{d}_s + \mathbf{A}\boldsymbol{\xi}_s - \mathbf{B} z_0)^{\mathrm{T}} \boldsymbol{u}_{s,\mathrm{fea}}^* \right)^2 \tag{6d}$$

where $\mathrm{N}_3$+$\mathrm{M}_3$ qubits $z_{3,i}$ and $\mathrm{N}_4$+$\mathrm{M}_4$ qubits $z_{4,ku}$ with the weight of $2^r$ and $2^u$ relax aggregated optimal cut (4b) and $S_{\mathrm{fea}}$ feasible cuts (4c) from inequality form into equality form, respectively; $\boldsymbol{\xi}_s=[\xi(z_1)]_{h\times 1}\leftarrow$UQAE.

## B. *Quantum-Based Feasible Cut Selection*

Although the optimal cuts can be aggregated as (4b) without breaking optimality, the feasible cuts in (4c) cannot because the aggregated feasible cut will no longer be able to exclude all identified infeasible trial solutions of ILP (4). As a result, Hamiltonian $H_3$ in (6d) that includes $S_{\mathrm{fea}}(\mathrm{N}_4+\mathrm{M}_4)$ superimposed qubits $z_{4,ku}$ will be very high-dimensional, which is unacceptable in the noisy intermediate-scale quantum era. To reduce the number of feasible cuts, we formulate a minimum set cover problem (7) to only select effective feasible cuts.

$$\min_{\boldsymbol{\sigma}_0} \mathbf{1}^{\mathrm{T}} \boldsymbol{\sigma}_0 \tag{7a}$$

$$\text{s.t. } \mathbf{N}\boldsymbol{\sigma}_0 \geq \mathbf{1} \tag{7b}$$

$$\boldsymbol{\sigma}_0 \in \{0,1\}^{S_{\mathrm{fea}}} \tag{7c}$$

where $\mathbf{N}$ is a matrix with $q\times S_{\mathrm{fea}}$ dimensions; $q$ is the cumulative number of identified infeasible trial solutions from ILP (4); $\boldsymbol{\sigma}_0$ represents a $S_{\mathrm{fea}}$ dimensional binary variable whose element equals 1 indicates that the corresponding feasible cut is selected, otherwise it is not selected.

Obviously, problem (7) is also a NP-hard ILP for CPU. Thus, similar to the Hamiltonian $H_{\mathrm{MP}}$ for ILP (4), we transform ILP (7) into quantum-amenable QUBO (8) by Hamiltonian $H_{\mathrm{CS}}$, which is composed of two Hamiltonians $H_4$ and $H_5$ that corresponds objective (7a) and constraint (7b), respectively.

$$\min_{\boldsymbol{\sigma}_0, \sigma_{1,rv}} H_{\mathrm{CS}} = H_4 + H_5 \tag{8a}$$

$$H_4 = \mathbf{1}^{\mathrm{T}} \boldsymbol{\sigma}_0 \tag{8b}$$

$$H_5 = \lambda_3 \sum_{r=1}^{q} \left( \mathbf{N}_r \boldsymbol{\sigma}_0 - 1 - \sum_{v=-\mathrm{M}_5}^{\mathrm{N}_5} 2^v \sigma_{1,rv} \right)^2 \tag{8c}$$

where $N_5+M_5$ qubits $\sigma_{1,rv}$ with the weight of $2^v$ relax inequality (7b) into equality form; $\lambda_3$ is an appropriate penalty factor.

*Remark 3:* Existing methods (e.g., [10]) select cuts fed back from one subproblem corresponding to multiple suboptimal trial solutions of the mater problem. Different from them, the proposed quantum-based feasible cut selection aims at removing redundant cuts from massive subproblems that correspond to one optimal trial solution of the mater problem.

### C. *Quantum-Assisted Iterative Process*

Fig. 3 illustrates the quantum-assisted iterative solution process for the two-stage stochastic optimization, where QPU and CPU are coordinated for computation so as to fully leverage their respective advantages.

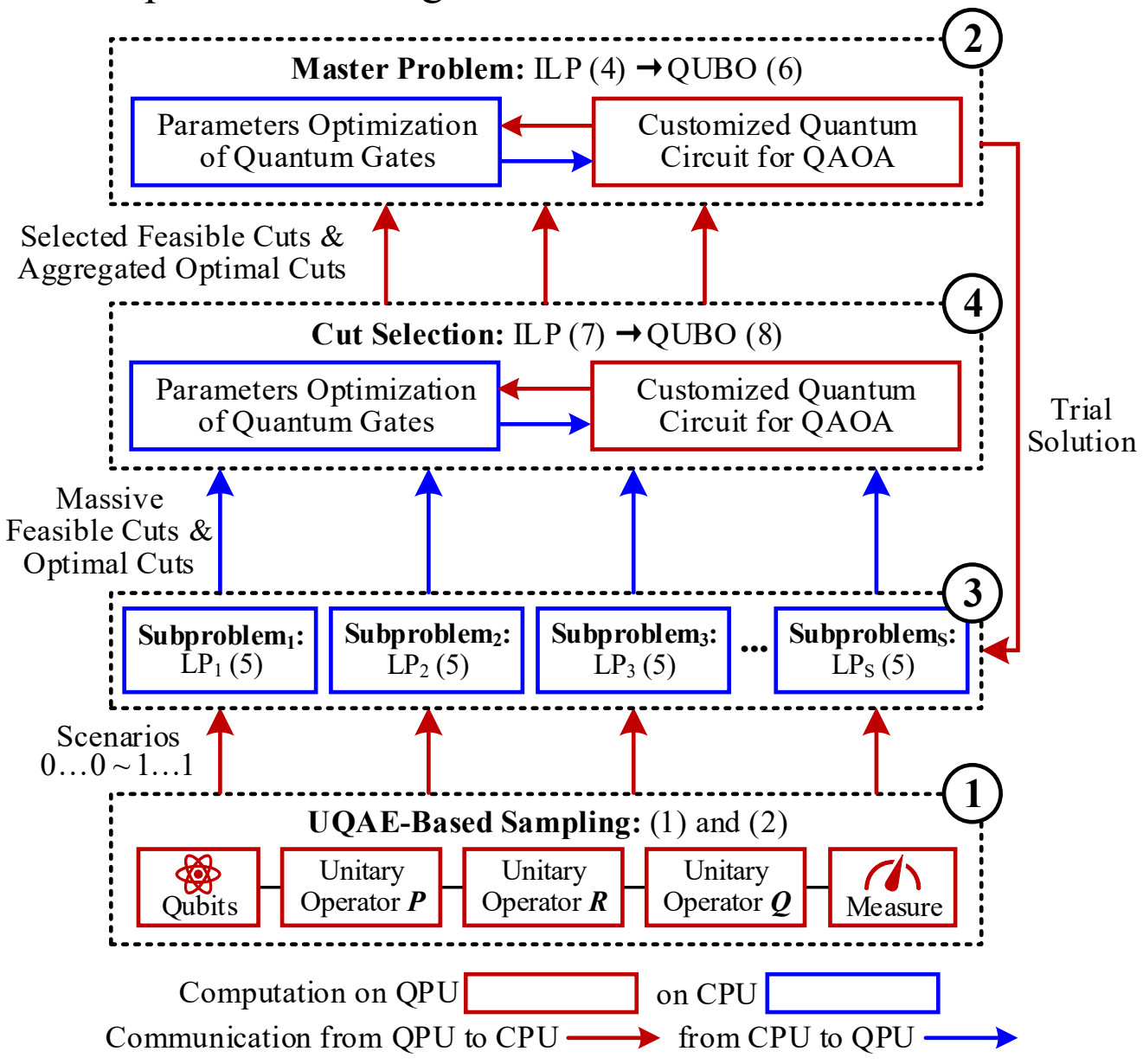

Fig. 3. Quantum-assisted iteration for two-stage stochastic optimization.

The detailed iterations are stated as follows:

*Step 1:* Generate massive scenarios on QPU. Model the distributions of RES prediction errors by unitary operator $\boldsymbol{P}$ in (1) and then implement UQAE-based sampling by unitary operators $\boldsymbol{R}$ and $\boldsymbol{Q}$ in (2) to generate S uncertain scenarios for subproblems. Go to Step 2.

*Step 2:* Solve the master problem on QPU+CPU. Add the selected feasible cuts and aggregated optimal cuts (if any) to ILP (4), which are further reformulated into QUBO (6). Then, search for the lowest energy of $H_{\mathrm{MP}}$ by QAOA with the customized quantum circuit to generate the trial solution. Go to Step 3.

*Step 3:* Solve the subproblems on CPU. Receive the trial solution from Step 1, and then solve independent LP (5) on CPU under S scenarios in parallel to feedback massive optimal and feasible cuts. Go to Step 4.

*Step 4:* Select cuts on QPU+CPU. Reformulate ILP (7) into QUBO (8) and thus search for the lowest energy of $H_{\mathrm{CS}}$ by QAOA. If there exists any new selected feasible cuts or aggregated optimal cuts at the current iteration, go back to Step 2 and start a new iteration; otherwise, terminate the iteration.

## IV. Numerical Experiments

Numerical experiments are conducted on Python 3.12 which employs IBM's Qiskit to design quantum circuits for UQAE and QAOA. IEEE 6-bus system is modified as a renewables-rich system to demonstrate the designed scheme, where three RESs are grid-connected at buses 1, 4, and 6, respectively; an energy storage system is equipped at bus 4; the flexibility of the two conventional generators at buses 1 and 6 has been modified. The dispatch horizon consists of 4 hours, in which each dispatch interval lasts 1 hour. The costs of start-up and generation for two conventional generators and the charging and discharging cost of energy storage system are considered under the constraints of direct current power flow, ramp limits, state of charge, etc.

Based on the quantum circuit for UQAE in Figs. 1 and 2, we only employ 3 qubits to model the distributions of RES prediction errors, and 1 ancilla qubit to sample $2^9$=512 times. The sample results are shown in Fig. 4. To achieve the same accuracy, sequential Monte Carlo simulation on CPU needs to sample $2^{10}$=1024 times. Thus, UQAE can lessen computational burdens for sampling with very few qubits and reduce sampling times. After the Cartesian product of discretized RES prediction errors, $2^3\times2^3\times2^3$=512 scenarios are generated with different weights (frequencies) for the stochastic optimization.

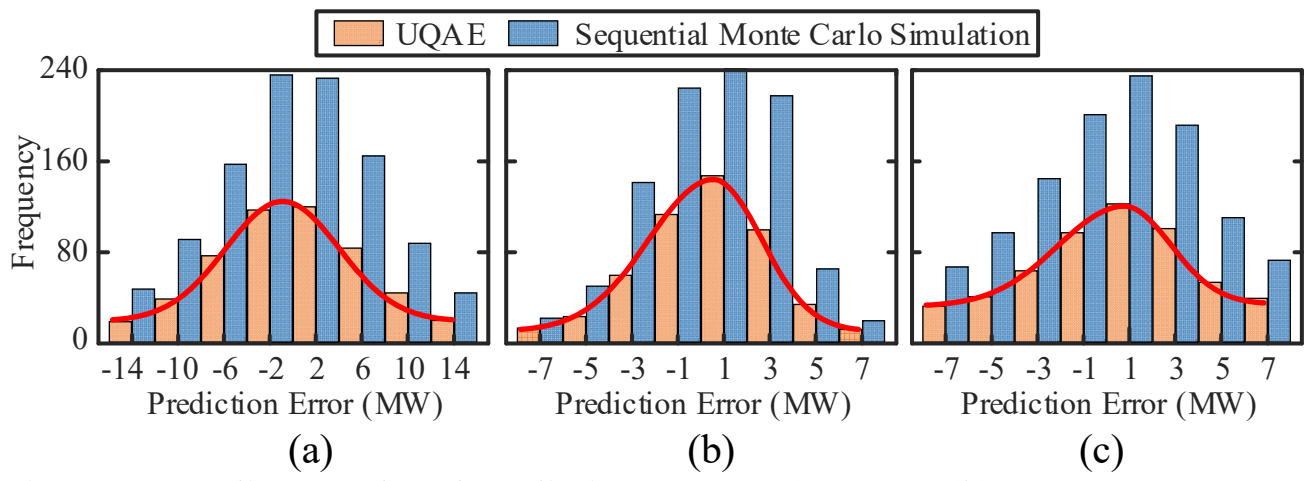

Fig. 4. Sampling results of prediction errors. (a) RES 1; (b) RES 2; (c) RES 3.

After integrating the scenarios, we compare the solution process under the quantum-assisted scheme and conventional CPU-based scheme. As shown in Fig. 5, based on quantum-assisted multi-cut Benders decomposition, both the optimal value and the iterative convergence process of upper/lower bound on QPU+CPU are almost the same as these on CPU. Notably, the difference in the minimum operation costs of 1770.76$ and 1770.79$ in Figs. 5(a) and 5(b) does not indicate that the quantum-assisted scheme has a loss of accuracy than the CPU-based one. Since the weight for each scenario has a certain randomness, according to the sample results in Figs. 4(a)-4(c), whether based on UQAE or sequential Monte Carlo simulation.

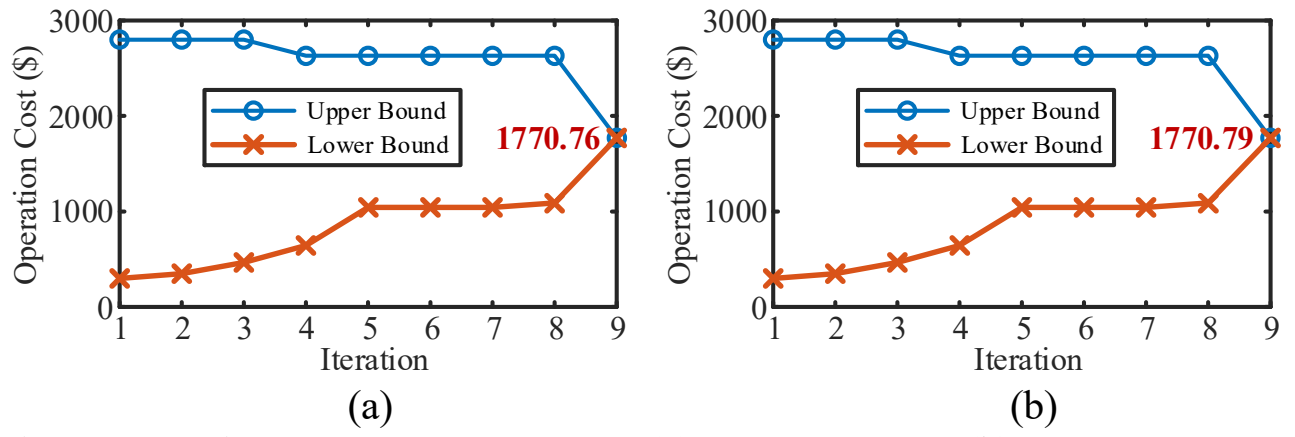

Fig. 5. Iterative convergence process. (a) on QPU+CPU; (b) on CPU.

Also, the quantum-assisted scheme will not change the economic dispatch results by stochastic optimization under all scenarios (see Fig. 6). Note that the negative values in Fig. 6(c) indicate the charging of energy storage system. As shown in Figs. 6(a)-6(c), the conventional generators and energy storage system balance the uncertain power flow in coordination with time-varying power outputs. Under some specific scenarios, the power output of conventional generator 1 ramps to the limit at 1h to cover the output gap of RESs. Also, the power output of conventional generator 2 lows to the limit or even shuts down to save the cost at 4h. Besides, the energy storage system shifts RES generations from 1h to 2h and 3h, thereby suppressing RES fluctuations and avoiding possible power congestion.

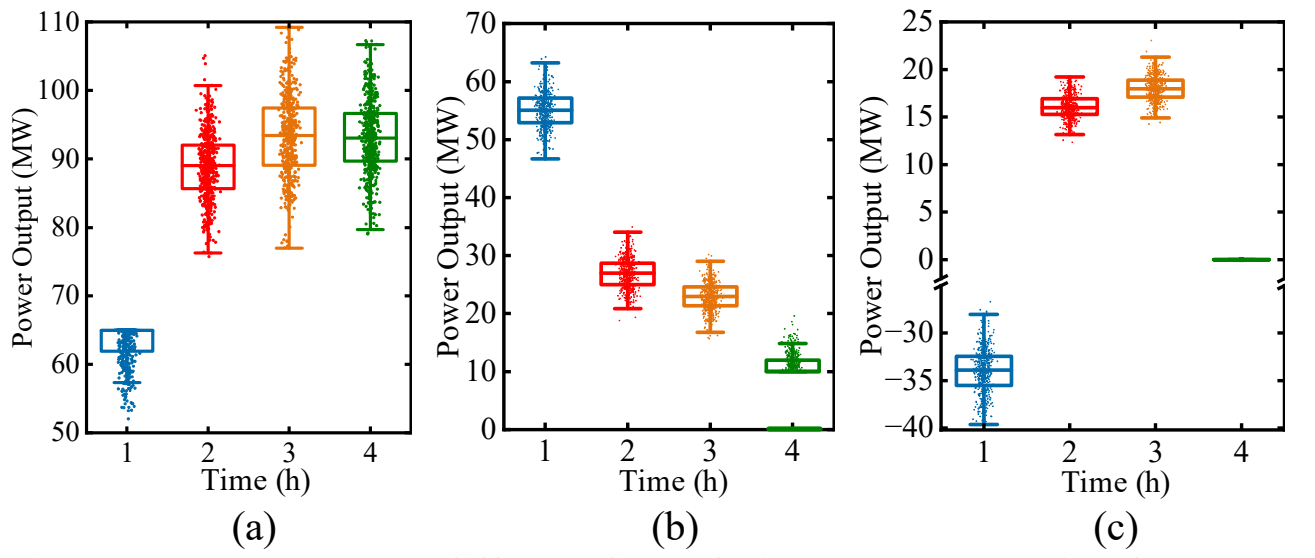

Fig. 6. Power outputs at different dispatch times. (a) Conventional generator 1; (b) Conventional generator 2; (c) Energy storage system.

In addition to the above correctness verification, Table I manifests the computational performance of multi-cut Bender decomposition under the proposed quantum-assisted scheme. In particular, UQAE can generate massive scenarios in only 0.11s. Although CPU-based commercial solvers (e.g., Gurobi) support multi-threaded parallel computing, the number of parallels is usually much smaller than the number of subproblems, limited by CPU cores. The total solution time of LP subproblems still rises to 16.32s. Benefiting from the fast optimization process for the customized Hamiltonians in QAOA, the total solution time of ILP master problem is controlled to the equivalent level of LP subproblems. Also, the quantum-based feasible cut selection is lightweight and effective, which only takes 2.23s to solve the NP-hard problem and thus reduces a total of 782 feasible cuts to 16 in the iterations. Hence, we can conclude that the quantum-assisted scheme is suitable for accelerating stochastic optimization. For larger-scale systems, the time complexities of these NP-hard problems will exponentially increase on the CPU, while the time complexities based on the quantum-assisted scheme are expected to increase in polynomials. Also, the quantum-assisted scheme should be modified to be parallel to reduce the required qubits and gates.

TABLE I. COMPUTATIONAL PERFORMANCE

| **Computing Time (s)** | | | |
|---|---|---|---|
| ***UQAE to Generate Scenarios*** | ***QAOA for Master Problem*** | ***QAOA for Cut Selection*** | ***Gurobi to solve Subproblems*** |
| 0.11 | 15.35 | 2.23 | 16.32 |
| **Total Number of Benders Cuts in Iterations** | | | |
| ***Original Feasible Cuts*** | ***Selected Feasible Cuts*** | ***Original Optimal Cuts*** | ***Aggregated Optimal Cuts*** |
| 782 | 16 | 3314 | 7 |

Without loss of generality, we test the robustness of the designed scheme in a total of 500 independent numerical experiments. As shown in Fig. 7, under either the quantum-assisted scheme on QPU+CPU or the conventional scheme on CPU, the operation cost fluctuates with the experiments. Since the weights (frequencies) of different scenarios driven from the random sample process of UQAE or sequential Monte Carlo simulation are some differences in each experiment. Fortunately, UQAE is more robust by the uniform discretization of unitary operator $\boldsymbol{R}$ and thus helps to simulate RES prediction errors accurately. Also, the average operation cost does not deviate from the optimal value. Hence, the quantum-assisted scheme performs better robustness than the CPU-based one.

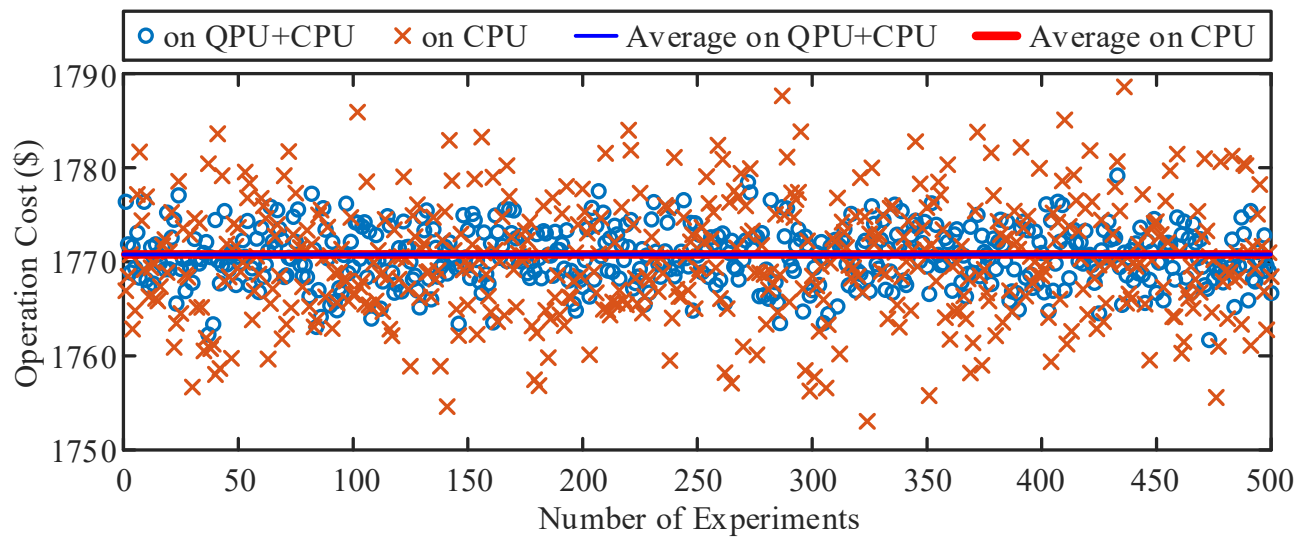

Fig. 7. Robustness test of the quantum-assisted scheme.

## V. Conclusions

This paper provides a quantum-assisted scheme to achieve fast, accurate, and robust stochastic optimization on the economic dispatch of renewables-rich power systems under massive scenarios. Leveraging the specific advantages of UQAE and QAOA, time-consuming stochastic optimization on CPU, including massive scenario generation and pertinent problem solution, is accelerated on QPU+CPU without losing accuracy. The quantum-assisted scheme paves the way for the uncertain optimization of power systems within polynomial complexity. Constrained by the limited number of qubits in the noisy intermediate-scale quantum era, the extension of this scheme in parallel for large-scale power systems is left for future studies.